\documentclass[a4paper,11pt]{amsart}
\usepackage{amsmath,amsthm,amssymb}
\usepackage[mathscr]{eucal}
 \usepackage{cite}
\usepackage{upgreek}
\usepackage[bookmarks,bookmarksnumbered]{hyperref}
\usepackage{enumerate}
\usepackage{amsmath}
\usepackage{mathrsfs}
\usepackage{blindtext}
\usepackage{scrextend}
\usepackage{enumitem}
\addtokomafont{labelinglabel}{\sffamily}
\usepackage{color}

\numberwithin{equation}{section}

\theoremstyle{plain}
\newtheorem{main}{Theorem}
\newtheorem{mcor}[main]{Corollary}

\newtheorem{theorem}{Theorem}[section]

\newtheorem{lemma}[theorem]{Lemma}

\theoremstyle{definition}

\newtheorem*{definition*}{Definition}

\newtheorem{remark}[theorem]{Remark}

\newcommand{\R}{\mathbb{R}}

\newcommand{\bS}{\mathbb{S}}

\newcommand{\cC}{\mathcal{C}}

\newcommand{\cM}{\mathcal{M}}

\newcommand{\cU}{\mathcal{U}}

\newcommand{\fg}{\mathfrak{g}}

\renewcommand{\d}{\text{d}}
\newcommand{\tcF}{\tilde{\mathscr{F}}}

\newcommand{\SO}{\operatorname{SO}}

\newcommand{\Isom}{\operatorname{Isom}}

\newcommand{\actson}{\curvearrowright}
\newcommand{\eps}{\varepsilon}

\begin{document}

\title[Local spectral gap in the group of euclidean isometries]
{Local spectral gap in the group of euclidean isometries}

\author[R. Boutonnet]{R\'{e}mi Boutonnet}
\address{Institut de Math\'{e}matiques de Bordeaux, Universit\'e de Bordeaux, 351 cours de la Lib\'eration, 33405 Talence Cedex,  France}
\email{remi.boutonnet@math.u-bordeaux.fr}

\author[A. Ioana]{Adrian Ioana}
\address{Department of Mathematics, University of California San Diego, 9500 Gilman Drive, La Jolla, CA 92093, USA, and IMAR, Bucharest, Romania}
\email{aioana@ucsd.edu}
\thanks{A.I. was supported in part by NSF Career Grant DMS \#1253402 and a Sloan Foundation Fellowship.}
\begin{abstract} 
We provide new examples of translation actions on locally compact groups with the ``local spectral gap property'' introduced in \cite{BISG15}.
 This property has applications to  strong ergodicity, the Banach-Ruziewicz problem,  orbit equivalence rigidity, and equidecomposable sets.
The main group of study here is the group $\Isom(\R^d)$ of orientation-preserving isometries of the euclidean space $\R^d$, for $d \geq 3$. 
We prove that the translation action of a countable dense subgroup $\Gamma$ on Isom$(\mathbb R^d)$ has local spectral gap, whenever the translation action of the rotation projection of $\Gamma$ on $\text{SO}(d)$ has spectral gap.
Our proof relies on the amenability of $\Isom(\R^d)$ and on work of Lindenstrauss and Varj\'u, \cite{LV14}.
\end{abstract}

\maketitle

\section{Introduction}
The main goal of this paper is to establish the so-called {\it local spectral gap} property for a new class of left translation actions $\Gamma \actson (G,m_G)$, where $\Gamma$ is a countable dense subgroup of a locally compact group $G$, any $m_G$ is a fixed left-invariant Haar measure of $G$. This notion was introduced by Salehi-Golsefidy and the authors in \cite{BISG15}. Before giving its precise definition and explaining some of its applications, let us discuss the classical notion of spectral gap for probability measure preserving actions.

A measure preserving action $\Gamma\curvearrowright (X,\mu)$ of a countable group $\Gamma$ on a standard probability space $(X,\mu)$ is said to have {\it spectral gap} if there exists a finite set $F\subset\Gamma$ such that the operator $T$ on $L^2(X)$ given by $T\varphi=\frac{1}{|F|}\sum_{g\in F}g\cdot\varphi$
satisfies $\|T\|_{L^2_0(X)}<1$. Here, we denote by $g\cdot\varphi$ the function $g\cdot\varphi(x)=\varphi(g^{-1}x)$ and by $L^2_0(X)\subset L^2(X)$ the subspace of functions orthogonal to the constants. Note that the condition $\|T\|_{L^2_0(X)}<1$  is equivalent to the existence of a constant $\kappa>0$ such that
\begin{equation}\tag{SG}
\text{$\|\varphi\|_2\leq\kappa\max_{g\in F}\|g\cdot\varphi-\varphi\|_2$, for every $\varphi\in L^2_0(X)$.}
\end{equation}
Spectral gap for probability measure preserving actions is an important tool with a wide range of applications to various areas of mathematics. 
Given a compact group $G$,  its Haar measure $m_G$ can be taken a probability measure, and thus the left translation action $\Gamma \actson G$ of any countable dense subgroup $\Gamma$ is probability measure preserving. The spectral gap property has been studied intensively for such translation actions, especially in the case where $G$ is a Lie group (see the introduction of \cite{BISG15} for a detailed account). 
The first results in this direction were obtained in the 1980s in connection with the Banach-Ruziewicz problem.
Thus, it was shown that for every $d\geq 3$ there exists a countable dense subgroup $\Gamma<\text{SO}(d)$ such that the left translation action $\Gamma\curvearrowright \text{SO}(d)$ has spectral gap \cite{Ma80,Su81,Dr84}.
In 2006, a breakthrough was made by Bourgain and Gamburd who proved that the left translation action $\Gamma\curvearrowright \text{SO}(3)$ has spectral gap, whenever $\Gamma<\text{SO}(3)$ is a dense subgroup generated by matrices with algebraic entries \cite{BG06}. 
This result has been subsequently generalized in \cite{BG10} to $\text{SU}(d)$, for every $d\geq 2$.  Most recently, it has been extended by Benoist and de Saxc\'{e} to cover arbitrary connected compact simple Lie groups 
\cite{BdS14}.

In the case when the group $G$  is non-compact,  the Haar measure $m_G$ is still invariant under translation actions of subgroups of $G$, but it is no longer finite. Therefore, the notion of spectral gap does not formally make sense in this setting. As we explain in \cite{BISG15}, a good extension is the following {\it local} notion, obtained by restricting attention to a fixed subset of $G$ with finite measure.

\begin{definition*}[\!\!{\cite{BISG15}}]
Consider an arbitrary standard measure space $(X,\mu)$ together with a measurable subset $B \subset X$ such that $0 < \mu(B) < \infty$. A measure preserving action $\Gamma \curvearrowright (X,\mu)$ is said to have {\it local spectral gap} with respect to $B$ if there exist a finite set $F\subset\Gamma$ and a constant $\kappa>0$ such that
\begin{equation}
\tag{LSG}
\text{$\|\varphi\|_{2,B}\leq\kappa\max_{g\in F}\|g\cdot\varphi-\varphi\|_{2,B}$, for every $\varphi\in L^2(X)$ with $\int_B\varphi\;\text{d}\mu=0$.}\end{equation}
Here, $\displaystyle{\|\varphi\|_{2,B}:=\Big(\int_B|\varphi|^2\;\text{d}\mu\Big)^{1/2}}$ denotes the $L^2$-norm of the restriction of $\varphi$ to $B$.
\end{definition*}

Here we say that a measure space $(X,\mu)$ is  {\it standard} if $X$ is a standard Borel space (i.e. a Polish space endowed with its Borel $\sigma$-algebra) and $\mu$ is a $\sigma$-finite Borel measure on $X$.

The main result of \cite{BISG15} is a generalization  of the above mentioned results from \cite{BG06,BG10,BdS14}  to the non-compact setting. Specifically, let $G$ be a connected simple Lie group, and $B\subset G$ be a bounded measurable set with  non-empty interior. It is proved in \cite{BISG15} that the left translation action $\Gamma\curvearrowright G$ has local spectral gap with respect to  $B$, whenever $\Gamma<G$ is a dense subgroup whose image through the adjoint representation of $G$ consists of matrices with algebraic entries (after fixing a basis of the Lie algebra $\fg$ of $G$).

In this paper we study the local spectral gap property for countable dense subgroups of the group $\text{Isom}(\mathbb R^d)=\mathbb R^d\rtimes\text{SO}(d)$ of orientation-preserving isometries of $\mathbb R^d$. If $d\in\{1,2\}$, then  $\text{Isom}(\mathbb R^d)$ is nilpotent and thus the local spectral gap property never holds. If $d\geq 3$, then one approach would be to follow the lines of \cite{BISG15} in order to establish the local spectral gap property for any dense subgroup of $\text{Isom}(\mathbb R^d)$ generated by isometries whose rotation and translation components have algebraic entries.
However, in this paper, we follow a different approach inspired by and building on the recent work of Lindenstrauss and Varj\'{u} \cite{LV14}, and prove the following stronger result.  In fact, our proof also gives a simpler approach to the main technical result of \cite{LV14} (see Section \ref{3}).

\begin{main}\label{A}
Let $G = \Isom(\R^d)$, for $d \geq 3$, and consider a countable dense subgroup $\Gamma < G$. Assume that the left translation action $\theta(\Gamma)\curvearrowright\SO(d)$ has spectral gap, where $\theta:\Isom(\mathbb R^d)\rightarrow\SO(d)$ denotes the natural quotient map. Let $\sigma$ denote either the left translation action $\Gamma\curvearrowright G$ or the natural isometric action $\Gamma \actson \R^d$.

Then $\sigma$ has local spectral gap with respect to any bounded measurable set $B$ with non-empty interior. 
\end{main}

In combination with results from \cite{BG06,BdS14}, Theorem \ref{A} implies that the actions $\Gamma\curvearrowright G$ and $\Gamma \actson \R^d$ have local spectral gap, for any countable dense subgroup $\Gamma<G$ whose rotation projection $\theta(\Gamma)$ consists of matrices with algebraic entries.

As one can check, the local spectral gap property for the quotient action $\Gamma \actson \R^d$ follows from that of the action $\Gamma \actson G$.
The local spectral gap property for the latter action, and thus Theorem \ref{A}, will be deduced from the following general result. 

\begin{main}\label{B} 
Consider $\Gamma < G$ as in Theorem \ref{A}. Denote by $\lambda:G\rightarrow\mathcal U(L^2(G))$ the left regular representation of $G$ given by $\lambda(g)\varphi(x)=\varphi(g^{-1}x)$, for every $g,x\in G$ and $\varphi\in L^2(G)$.
 
 Then for any compact set $K\subset G$, we can find  a finite set $F\subset\Gamma$ and a constant $\kappa>0$ such that 
 \[\text{$\sup_{g\in K}\|\lambda(g)\varphi-\varphi\|_2\leq \kappa\max_{g\in F}\|\lambda(g)\varphi-\varphi\|_2$, for every $\varphi\in L^2(G)$}.\]
\end{main}

Let us put this result into perspective. 
Recall that the amenability of a locally compact group $G$ can be characterized by the existence of a net $\{\varphi_n\}_n$ of unit vectors in $L^2(G)$ satisfying either one of the following three ``almost invariance" conditions (see e.g. \cite[Appendix G]{BdHV08}):
\begin{enumerate}[label=(\roman*)]
\item ${\sup_{g\in K}\|\lambda(g)\varphi_n-\varphi_n\|_2\rightarrow 0}$, for every compact subset $K\subset G$.
\item $\|\lambda(g)\varphi_n-\varphi_n\|_2\rightarrow 0$, for every $g\in G$.
\item $\|\lambda(g)\varphi_n-\varphi_n\|_2\rightarrow 0$, for every $g$ belonging to a dense subset $S\subset G$.
\end{enumerate}
In particular, the existence of a net satisfying (ii) or  (iii) implies the existence of a net satisfying (i).  However, as can be easily seen in the case $G=\mathbb R^k$, it is not true in general that any net verifying  (ii) or (iii) must verify (i). 
Nevertheless, Theorem \ref{B} implies that in the case $G=\text{Isom}(\mathbb R^d)$, $d\geq 3$, any net of unit vectors $\{\varphi_n\}_n$ in $L^2(G)$ satisfying (ii) or (iii) for $S=\Gamma$ must necessarily satisfy (i). 

One can now deduce Theorem \ref{A} by making the key observation that, since $G$ is amenable, the last fact implies the local spectral  gap property for $\Gamma\curvearrowright G$. This observation is inspired by Margulis' work \cite{Ma82}, and is made precise in Lemma \ref{BtoA}.

Similarly, one can deduce the local spectral gap property for $\Gamma\curvearrowright\mathbb R^d$ whenever the estimate from Theorem \ref{B} holds for the natural representation $\lambda_0:G\rightarrow\mathcal U(L^2(\mathbb R^d))$ in place of $\lambda$. Note that in the setting of Theorem \ref{A} this estimate and thus the local spectral gap property for $\Gamma\curvearrowright\mathbb R^d$ can be alternatively derived by using the main technical result of \cite{LV14}.

Next, we record several general consequences of local spectral gap observed in our particular setting. 
\begin{mcor}\label{C}
Consider $\Gamma < G$ as in Theorem \ref{A}. Then the following hold true.
\begin{enumerate}
\item The left translation action $\Gamma\curvearrowright G$ is strongly ergodic.
\item Denote by $\cC(G)$ the family of measurable subsets $A\subset G$ with compact closure. Then the Haar measures of $G$ are the only finitely additive measures $\nu: \cC(G) \to [0,\infty)$ that are $\Gamma$-invariant.
\item Let $\Lambda$ be any countable dense subgroup of a connected Lie group $H$ with trivial center. Then the left translation actions $\Gamma \actson G$ and $\Lambda \actson H$ are orbit equivalent if and only if there exists a topological isomorphism $\delta:G \to H$ such that $\delta(\Gamma) = \Lambda$.
\item Consider two measurable subsets $A, B \subset G$ with non-empty interior and compact closure. Then $A$ and $B$ are measurably $\Gamma$-equidecomposable if and only if they have the same Haar measure.
\end{enumerate}
Moreover, the analogous conclusions to (1), (2) and (4) for the action $\Gamma \actson \R^d$ also hold.
\end{mcor}

Corollary \ref{C} is obtained by combining Theorem \ref{A} with results from \cite{BISG15} (parts (1) and (2)), \cite{Io14} (part (3)), and \cite{GMP16} (part (4)). 
Towards recalling the notions involved in its statement, consider a measure preserving action $\Gamma \actson (X,\mu)$ of a countable group on a standard measure space.

 In order to recall what it means for the action to be strongly ergodic, since the measure $\mu$ can be infinite, we first choose a probability measure $\mu_0$ on $X$ with the same null sets as $\mu$. 
The action is said to be {\it strongly ergodic} 
if any sequence $\{A_n\}_n$ of measurable subsets of $X$ satisfying $\mu_0(g\cdot A_n\;\Delta\; A_n)\rightarrow 0$, for all $g\in\Gamma$, must satisfy $\mu_0(A_n)(1-\mu_0(A_n))\rightarrow 0$ \cite{CW80,Sc80}. It is easy to see that this definition does not depend on the choice of $\mu_0$.

We also recall that $\Gamma \actson (X,\mu)$ is {\it orbit equivalent} to another measure preserving action $\Lambda \actson (Y,\nu)$  if there exists a measure class preserving Borel isomorphism $\theta: X \to Y$ such that $\theta(\Gamma \cdot x) = \Lambda \cdot \theta(x)$, for $\mu$-almost every $x \in X$. Thus, (2) can be interpreted to mean that one can completely recognize the translation action $\Gamma \actson G$, up to conjugacy, from the measurable structure of its orbits.

Finally, recall that $A$ and $B$ are $\Gamma$-{\it equidecomposable} if and only if there exist finite partitions into subsets $A = \sqcup_{i=1}^n A_i$ and $B = \sqcup_{i = 1}^n B_i$ such that $B_i = g_i \cdot A_i$ for some $g_i \in \Gamma$, for every $i = 1,\dots,n$. They are said to be {\it measurably} equidecomposable if each piece in the above partitions is measurable.


\subsection*{Acknowledgement}

Part of this work was done while the first author was visiting the University of California at San Diego. He thanks the Mathematics Department for the kind hospitality he received. The authors are very grateful to the anonymous referee for several helpful comments, and in particular for asking the question which led to Lemma \ref{converse} and for pointing out a gap in our original proof of Lemma \ref{equideco}.

\section{Preliminaries}

\subsection{General notations}

Given any locally compact group $H$, we denote by $L^2(H)$ the Hilbert space of square integrable functions on $H$ with respect to the Haar measure. 

We denote by $\mathcal M(H)$ the family of Borel probability measures on $H$, and by $C_c(H)$ the space of compactly supported continuous functions on $H$.
If $\mu,\nu\in\mathcal M(H)$, their convolution product $\mu*\nu\in\mathcal M(H)$ is defined by  
\[\int_{H}f\;\text{d}(\mu*\nu) = \int_{H}\int_{H}f(xy)\;\text{d}\mu(x)\text{d}\nu(y), \text{ for every } f\in C_c(H).\]
For $n\geqslant 1$,  $\mu^{*(n)}$ denotes the $n$-fold convolution product of $\mu$ with itself. 
We denote by $\check{\mu}$ the probability measure given by 
\[\int_{H}f\;\text{d}\check{\mu}=\int_{H}f(x^{-1})\;\text{d}\mu(x), \text{ for every } f\in C_c(H),\] and say that $\mu$ is {\it symmetric} if  $\check{\mu}=\mu$.

Given a unitary representation $\pi:H\rightarrow\mathcal U(\mathcal H)$ and $\mu\in\mathcal M(H)$, the formula $\pi(\mu)=\int_{H}\pi(g)\;\text{d}\mu(g)$ defines a bounded linear operator on $\mathcal H$ with $\|\pi(\mu)\|\leq 1$. Then $\pi(\mu)^*=\pi(\check{\mu})$, and thus $\pi(\mu)$ is self-adjoint whenever $\mu$ is symmetric. Moreover, $\pi(\mu*\nu)=\pi(\mu)\pi(\nu)$, for every $\nu\in\mathcal M(H)$.

\subsection{A direct integral decomposition}

Let $d\geq 3$ and denote by $G=\text{Isom}(\mathbb  R^d)=\mathbb R^d\rtimes\text{SO}(d)$ the group of orientation-preserving isometries of $\mathbb R^d$.  Thus, every isometry $g\in G$ can be uniquely written as $g=(v(g),\theta(g))$, where $v(g)\in\mathbb R^d$ and $\theta(g)\in \SO(d)$. For $g\in G$ and $x\in\mathbb R^d$, we then have $g(x)=v(g)+\theta(g)x$.
 The product of two isometries is given by $(v_1,\theta_1)(v_2,\theta_2)=(v_1+\theta_1(v_2),\theta_1\theta_2)$. 

Consider the left regular representation $\lambda: G\rightarrow\mathcal U(L^2(G))$ of $G$ given by $\lambda(g)\varphi(x)=\varphi(g^{-1}x)$. 
Consider the unitary operator $U:L^2(\SO(d))\otimes L^2(\mathbb R^d)\rightarrow L^2(G)$ given by $(Uf)(\xi,\omega)=f(\omega,\omega^{-1}\xi)$, for every $f\in L^2(\SO(d))\otimes L^2(\mathbb R^d)$, $\omega\in\SO(d)$, and $\xi\in\mathbb R^d$. Let $g=(v,\theta)\in G$. We claim that \begin{equation}
\label{regular}\big(U^{-1}\lambda(g)U\big)f(\omega, \xi)=f(\theta^{-1}\omega,\xi-\omega^{-1}v),\end{equation}
for every $f\in L^2(\SO(d))\otimes L^2(\mathbb R^d)$, $\omega\in \SO(d)$, and $\xi\in\mathbb R^d$. To see this, let $\tilde f\in L^2(\SO(d))\otimes L^2(\mathbb R^d)$  be given by $\tilde f(\omega,\xi)=f(\theta^{-1}\omega,\xi-\omega^{-1}v)$. Then for every $\omega\in\SO(d)$ and $\xi\in\mathbb R^d$, we have that
\begin{align*}
\lambda(g)Uf(\xi,\omega)&=Uf(g^{-1}(\xi,\omega))=Uf((-\theta^{-1}v,\theta^{-1})(\xi,\omega))\\&=Uf(-\theta^{-1}v+\theta^{-1}\xi,\theta^{-1}\omega)=f(\theta^{-1}\omega,-\omega^{-1}v+\omega^{-1}\xi)\\&=\tilde f(\omega,\omega^{-1}\xi)=U\tilde f(\xi,\omega),
\end{align*}
which proves \eqref{regular}.

Now, we rewrite the representation $\lambda$ using the Fourier transform on $\mathbb R^d$.  
We let $e(y):=e^{-i2\pi y}$, for $y\in\mathbb R$, and denote by $\langle \cdot ,\cdot\rangle$ and $|x|=\langle x,x \rangle^{1/2}$, for $x\in\mathbb R^d$, the usual inner product and norm on $\mathbb R^d$. We  denote by $\mathscr F:L^2(\mathbb R^d)\rightarrow L^2(\mathbb R^d)$ the Fourier transform given by 
\[\mathscr F(f)(x)=\int f(\xi)e(\langle x,\xi\rangle)\;\text{d}\xi, \text{ for every } f\in L^2(\mathbb R^d) \text{ and } x\in\mathbb R^d.\]

We denote still by  $\mathscr F$ 
the operator $\text{Id}\otimes\mathscr F$ on $L^2(\SO(d))\otimes L^2(\mathbb R^d)$. For $g=(v,\theta)\in G$, we denote $\tilde\lambda(g)=U^{-1}\lambda(g)U\in\mathcal U(L^2(\SO(d))\otimes L^2(\mathbb R^d))$.
Let  $f\in L^2(\SO(d))\otimes L^2(\mathbb R^d)$. By using \eqref{regular} we get 
\begin{align*}
\mathscr F(\tilde\lambda(g)f)(\omega,x) 
&=\int f(\theta^{-1}\omega,\xi-\omega^{-1}v)e(\langle x,\xi\rangle)\;\text{d}\xi\\
&=e(\langle x,\omega^{-1}v\rangle)\int f(\theta^{-1}\omega,\xi)e(\langle x,\xi\rangle)\;\text{d}\xi\\
&=e(\langle \omega x,v\rangle)\mathscr F(f)(\theta^{-1}\omega,x).
\end{align*}
Therefore, we have $(\mathscr F\tilde\lambda(g)\mathscr F^{-1})f(\omega,x)=e(\langle \omega x,v\rangle)f(\theta^{-1}\omega,x)$, for every $f \in L^2(\SO(d))\otimes L^2(\mathbb R^d)$.
Thus, if we denote $\tcF=\mathscr F U^{-1}$, then we have $(\tcF\lambda(g)\tcF^{-1})f(\omega,x)=e(\langle \omega x,v\rangle)f(\theta^{-1}\omega,x)$, for every $f \in L^2(\SO(d))\otimes L^2(\mathbb R^d)$.

For $x\in\mathbb R^d$, we define a unitary representation $\pi_x:G\rightarrow\mathcal U(L^2(\SO(d))$ by letting 
\begin{equation}\label{pi}\pi_x(g)\varphi(\omega)=e(\langle \omega x,v\rangle)\varphi(\theta^{-1}\omega), \text{ for every } \varphi\in L^2(\SO(d)).\end{equation}

Thus, if we disintegrate $\displaystyle{L^2(\SO(d))\otimes L^2(\mathbb R^d)=\int_{\mathbb R^d}^{\oplus}L^2(\SO(d))\;\text{d}x}$, then \begin{equation}\label{disint}\tcF\lambda \tcF^{-1}=\int_{\mathbb R^d}^{\oplus}\pi_x\;\text{d}x.\end{equation}

\section{Main technical result}\label{3}

In this section we generalize an estimate due to Lindenstrauss and Varj\'u, \cite[Theorem 2.1]{LV14}. Our proof relies on ingredients already used in \cite{LV14} (see also \cite[Section 3]{CG11} for related ideas), but it is much simpler.  Specifically, we use certain elementary spectral gap estimates contained in \cite[Section 3]{LV14}, but neither flattening arguments nor Littlewood-Paley decompositions.

\begin{theorem}\label{maintech} Let $G=\Isom(\mathbb R^d)$, for $d\geq 3$, and
consider the representations $\pi_x$, $x \in \R^d$, introduced in the previous section.
Consider a symmetric probability measure $\mu\in\mathcal M(G)$ and define operators $T_x:=\pi_x(\mu)$ on $L^2(\SO(d))$, for all $x \in \R$. We make the following assumptions: 
\begin{enumerate}
\item There exists $\alpha > 0$ such that $\|T_0\varphi\|_2\leq (1 - \alpha) \|\varphi\|_2$, for every $\varphi \in L^2_0(\SO(d))$;
\item For every $x \in \R^d$, we have $\mu(\{g \in G \, \vert \, gx \neq x\}) > 0$.
\end{enumerate}
Then there exists a constant $c_0 > 0$ (depending on $d$ and $\mu$) such that 
\[\Vert T_x \Vert \leq 1 - c_0\min(\vert x \vert^2,1), \text{ for every } x \in \R^d.\]
\end{theorem}

\begin{remark}
Proceeding as in \cite{LV14}, we could express explicitly the constant $c_0$ in terms of $\mu$ under certain finite moment assumptions on the translation part of $\mu$. Since this aspect is not relevant for our purposes, we will not pursue it here.
\end{remark}

\begin{remark}\label{rem1}
Theorem 2.1 in \cite{LV14} proves an estimate of the form $\Vert S_r \Vert \leq 1 - c\min(r^2,1)$ for all $r \geq 0$, where $S_r := \rho_r(\mu)$ and $\rho_r$ denotes the representation of $G$ on $L^2(\bS^{d-1})$ defined by 
\[\rho_r(g)\varphi(\xi) = e(r\langle \xi,v(g)\rangle)\varphi(\theta(g)^{-1}\xi).\]
We point out that, for any $x \in \R^d$, the following inequality holds
\begin{equation}\label{StoT}\Vert S_{\vert x \vert} \Vert \leq \Vert T_x \Vert.\end{equation}
Indeed, denote by $K_x=\{\theta\in \SO(d)|\theta(x)=x\}$ the stabilizer of $x$ in $\SO(d)$ and by $L^2(\SO(d))^{K_x}\subset L^2(\SO(d))$ the subspace of right $K_x$-invariant functions. Since $\pi_x$ commutes with the right action of $K_x$, the subspace $L^2(\SO(d))^{K_x}$ is $\pi_x(G)$-invariant. Identifying $\SO(d)/K_x$ with $\bS^{d-1}$, we see that the restriction of $\pi_x$ to $L^2(\SO(d))^{K_x}$ is unitarily equivalent to $\rho_{|x|}$. Hence, \eqref{StoT} follows.
\end{remark}

For the rest of this section, we fix a symmetric probability measure $\mu$ on $G$. For $x \in \R^d$, we denote  by $T_x = \pi_x(\mu)$ the corresponding averaging operator on $L^2(\SO(d))$, where $\pi_x$ is defined in the previous section. We will bound $\Vert T_x \Vert$ by treating separately the case where $\vert x \vert$ is large and the case where  $\vert x \vert$ is small.

\begin{lemma}\label{large}
Assume that $\|T_0\varphi\|_2\leq\|\varphi\|_2/2$, for every $\varphi \in L^2_0(\SO(d))$.
Then there exists an absolute constant $C > 0$ such that for all $x,y \in \R^d$ satisfying $\vert x \vert \geq \vert y \vert/2$, we have
\[ 1 - \Vert T_y \Vert \leq C(1 - \Vert T_x \Vert).\]
\end{lemma}

{\it Proof.} We start with a claim.

{\bf Claim.} For all $a,b \in \R^d$ such that $\vert a \vert = \vert b \vert$, we have $\Vert T_a \Vert = \Vert T_b \Vert$.

Take an element $h \in \SO(d)$ such that $b = ha$. Denote by $\sigma: \SO(d) \to \cU(L^2(\SO(d)))$ the right regular representation, and observe that for all $\varphi \in L^2(\SO(d))$ and $\omega \in \SO(d)$,
\[\sigma(h)T_a\varphi(\omega) = (T_a\varphi)(\omega h)  = \int_G e(\langle \omega h a,v(g)\rangle) \varphi(\theta(g)^{-1}\omega h)\d \mu(g) = T_b\sigma(h) \varphi(\omega).\]
This computation gives that $T_b = \sigma(h)T_a\sigma(h)^{-1}$, and the claim follows.

Since $\vert x \vert \geq \vert y \vert/2$, we may find $x' \in \R^d$ such that $\vert x' \vert = \vert x \vert$ and $\vert x - x' \vert = \vert y \vert$. From the claim, we have $\Vert T_x \Vert = \Vert T_{x'} \Vert$ and $\Vert T_y \Vert = \Vert T_{x - x'} \Vert$.

Next, we use some of the arguments in \cite[Section 3]{LV14} to conclude. Given $\eps > 0$, we say that a function $\varphi \in L^2(\SO(d))$ with unit norm is {\it $\eps$-invariant} for $T_x$ if $\Vert T_x\varphi \Vert_2 \geq 1 - \eps^2$.

Let $\eps > 0$ such that $1-\eps^2 < \Vert T_x\Vert$. By the above claim we also have $1-\eps^2 < \Vert T_{x'}\Vert$. So we can find unit vectors $\varphi_1,\varphi_2\in L^2(\SO(d))$ which are $\eps$-invariant for $T_x$ and $T_{x'}$, respectively. 
Using the assumption on $T_0$ and arguing as in the proof of \cite[Lemma 3.3]{LV14} we find an absolute constant $C > 0$, and functions $\psi_1,\psi_2\in L^2(\SO(d))$ such that  $|\psi_i|\equiv 1$ and $\|\psi_i-\varphi_i\|_2\leq C\varepsilon$ for both $i = 1,2$, and $\psi_1$ and $\psi_2$ are $C\eps$-invariant for $T_x$ and $T_{x'}$, respectively. 

Observing that
\[T_{x-x'}(\psi_1\overline{\psi_2})=\int_G\pi_x(g)\psi_1\cdot\overline{\pi_{x'}(g)\psi_2}\;\text{d}\mu(g),\]
and repeating the proof of \cite[Lemma 3.4]{LV14}  implies that $\psi_1\overline{\psi}_2\in L^2(\SO(d))$ is a $C'\varepsilon$-invariant unit vector for $T_{x-x'}$, for some absolute constant $C'>0$. This leads to $\Vert T_{x - x'}\Vert \geq 1 - {C'}^2\eps^2$, hence
\[  1 - \Vert T_y \Vert = 1 - \Vert T_{x - x'} \Vert \leq {C'}^2\eps^2,\]
for any $\eps>0$ satisfying $1-\eps^2 < \Vert T_x\Vert$. This proves that $1-\|T_y\|\leq C'^2(1-\|T_x\|)$, as wanted.
\hfill $\blacksquare$



\begin{lemma}\label{small}
Assume that the measure $\mu$ is compactly supported, that $\int_G v(g)\d\mu(g) = 0$ and that $\mu(\{g\in G \, \vert \, v(g) \neq 0\}) > 0$. 
Assume moreover that  $\|T_0\varphi\|_2\leq\|\varphi\|_2/2$, for every $\varphi \in L^2_0(\SO(d))$.
Then there exists a constant $c > 0$ such that for all $x \in \R^d$ satisfying $|x| \leq c$, we have 
\[\Vert T_x \Vert \leq 1 - c|x|^2.\]
\end{lemma}

{\it Proof.} The proof follows the same lines as the proof of \cite[Lemma 3.6]{LV14}. Put $C:=\int_G|v(g)|^2\;\text{d}\mu(g)$. 
We claim that $\Vert T_x\varphi - T_0\varphi \Vert_2 \leq 2\pi C^{1/2}|x|\Vert \varphi \Vert_2$, for every $\varphi\in L^2(\SO(d))$. Indeed, we have
\begin{align*}
\Vert T_x\varphi - T_0\varphi \Vert_2^2 & = \int_{\SO(d)} \vert (T_x - T_0)\varphi(\omega)\vert^2\d \omega\\
& \leq \int_{\SO(d)}\int_G \vert e(\langle \omega x,v(g) \rangle) - 1 \vert^2\vert \varphi(\theta(g)^{-1}\omega)\vert^2 \d \mu(g) \d \omega\\
\end{align*}
Since $|e(\langle \omega x,v(g) \rangle) - 1 \vert^2 = 2 - 2\cos(2\pi\langle \omega x,v(g) \rangle) \leq 4\pi^2\vert v(g)\vert^2\vert x \vert^2$, for all $g\in G$ and $\omega\in\SO(d)$, the claim follows. By using the triangle inequality and the assumption on $T_0$, we get  that if $|x|$ is small enough, then
\begin{equation}\label{Ton0int}
\Vert T_x\varphi \Vert_2 \leq 2^{-1/2}\Vert \varphi \Vert_2, \text{ for all }\varphi \in L^2_0(\SO(d)).
\end{equation}

Now, since $|e(\langle\omega x,v(g)\rangle)-1+i2\pi\langle\omega x,v(g)\rangle|\leq 4\pi^2|v(g)|^2|x|^2$, for all $g\in G$ and $\omega\in\SO(d)$, by integrating we get that $\|T_x1-1\|_2\leq 4\pi^2C|x|^2$. Since $\|T_x\|\leq 1$, this further implies that \begin{equation}\label{sec}\|T_x^21-1\|_2\leq 2\|T_x1-1\|_2\leq 8\pi^2C|x|^2.\end{equation} 

Next, if we put $\tilde \mu := \mu \ast \mu$, then we have
\[\Vert T_x 1 \Vert_2^2 = \langle T_x^*T_x 1,1\rangle = \int_{\SO(d)}\int_G e(\langle \omega x,v(g)\rangle)\d\tilde \mu(g) \d \omega=  \int_{\SO(d)}\int_G \cos(2\pi\langle \omega x,v(g)\rangle)\d\tilde \mu(g)\d \omega,\]
where the last equality is obtained by taking real parts. 
The assumptions made on $\mu$ imply that 
 $\int_G \vert v(g) \vert^2 \d\tilde \mu(g)=2C$ (see \cite[Lemma 8.2]{LV14}).  By using that $\int_{\SO(d)} \langle \omega x,v\rangle^2 \d \omega = \vert x \vert^2\vert v \vert ^2/d$ and $\cos(a) \leq 1 - a^2/2 + a^4/24$, for all $x,v\in\mathbb R^d$ and $a \in \R$, and the last displayed formula, we get that
\begin{align*}\Vert T_x1 \Vert_2^2 - (1 - \frac{4\pi^2 C}{d}\vert x \vert^2) &=  \int_{\SO(d)}\int_G \cos(2\pi\langle \omega x,v(g)\rangle) - (1 - (2\pi\langle \omega x,v(g)\rangle)^2/2)\d \tilde \mu(g)\d \omega\\&\leq \frac{2\pi^4}{3}|x|^4 \int_G \vert v(g) \vert^4 \d\tilde \mu(g).\end{align*}
In particular, we deduce that if $|x|$ is small enough, then
\begin{equation}\label{third}
\|T_x1\|_2^2\leq 1-\frac{2\pi^2 C}{d}\vert x\vert^2.
\end{equation}
The conclusion now follows by combining  \eqref{Ton0int},  \eqref{sec}, \eqref{third}, and proceeding exactly as in the proof of \cite[Lemma 3.6]{LV14} with $T_x$ and $|x|$ instead of $S_r$ and $r$.
\hfill$\blacksquare$

We can now prove the theorem.

{\it Proof of Theorem \ref{maintech}.}
Let $\mu$ be a symmetric probability measure on $G$ satisfying assumptions (1) and (2). 
Following an idea from the proof of \cite[Lemma 9.1]{LV14}, we first reduce to the case when $\mu$ is compactly supported by restricting it to a sufficiently large ball.

For $s > 0$, we denote by $B_s \subset \R^d$ the ball of radius $s$ centered at $0$, and set $K_s := B_s\cdot \SO(d) \subset G$. Write $\mu_s^0$ for the restriction of $\mu$ to $K_s$ given by $\mu_s^0(A) = \mu(A \cap K_s)$ for all measurable sets $A \subset G$, and set $\nu_s^0 := \mu - \mu_s^0$. Finally, consider the normalized measures $\mu_s := \mu_s^0/\mu(K_s)$ and $\nu_s := \nu_s^0/\mu(G \setminus K_s)$.

Take $s > 0$ large enough so that $\beta:= \mu(G \setminus K_s) < \alpha/2$, where $\alpha$ is the constant appearing in assumption (1). Since $\alpha\leq 1$, we have $\mu(K_s)>0$. Since $\mu = (1 - \beta) \mu_s + \beta \nu_s$, by using (1) we get
\begin{equation}\label{gaps}
\Vert \pi_0(\mu_s)\varphi\Vert_2 \leq \frac{\Vert \pi_0(\mu)\varphi\Vert_2 + \beta\Vert \pi_0(\nu_s)\varphi\Vert_2 }{1 - \beta} \leq \frac{1 - \alpha + \beta}{1- \beta}\Vert \varphi \Vert_2, \text{ for all } \varphi \in L^2_0(\SO(d)).
\end{equation}
In particular, $\mu_s$ satisfies condition (1), possibly with a different choice of $\alpha > 0$.

We claim that there exists at most one point $x_0 \in \R^d$ such that $\mu_s(\{g \, \vert \, gx_0 \neq x_0\}) = 0$. Indeed, if we had two distinct such points $x_0$ and $x_1$, then $\mu_s$-almost surely, $\theta(g)$ would fix $x_0 - x_1$. In this case the $L^2$-function on $\SO(d)$ defined by $\omega \mapsto\langle\omega (x_0 - x_1), x_0 - x_1\rangle\in\mathbb C$ would be non-constant and $\pi_0(\mu_s)$-invariant, contradicting the spectral gap property \eqref{gaps}. This proves our claim.

If such an $x_0$ does not exist, then $\mu_s$ also satisfies condition (2).
If such an $x_0$ does exist, then since $\mu(\{g \, \vert \, gx_0 \neq x_0\}) > 0$, we may find $s' > s$ such that $\mu_{s'}(\{g \, \vert \, gx_0 \neq x_0\}) > 0$. 
This implies that $\mu_{s'}$ satisfies condition (2). Indeed, if $x_1\in \R^d$ is a point fixed by $\mu_{s'}$-almost every $g\in G$, then since $\mu_s\leq \kappa\mu_{s'}$ for some $\kappa>0$, we would get $x_1$ is fixed by $\mu_s$-almost every $g\in G$. By the previous paragraph, this would force that $x_1=x_0$, contradicting that $\mu_{s'}(\{g \, \vert \, gx_0 \neq x_0\}) > 0$.
Moreover, since $s' > s$ we have that $\mu(G\setminus K_{s'})<\alpha/2$ and the computation from \eqref{gaps} shows that $\mu_{s'}$ also satisfies condition (1).

In summary, we have found $s'>0$ such that $t:=\mu(K_{s'})>0$ and $\mu_{s'}$ satisfies conditions (1) and (2). 
Since $\mu = t\mu_{s'} + (1 - t)\nu_{s'}$, proving the conclusion for $\mu$ reduces to proving the conclusion for $\mu_{s'}$. 
Since $\mu_{s'}$ is compactly supported, it follows that we may assume that $\mu$ is compactly supported.

Next, since $T_0=\pi_0(\mu)$ satisfies the spectral gap condition (1), there is no non-zero vector $x \in \R^d$ such that $\theta(g)x = x$, for  every $g$ in the support of $\mu$. Indeed, otherwise $\SO(d)\ni\theta\mapsto\langle\theta x,x\rangle\in\mathbb C$ would be a non-constant $T_0$-invariant function belonging to $L^2(\SO(d))$.

By \cite[Lemma 4]{Va12}, there exists a unique $a \in \R^d$ such that $\int g(a)\;\text{d}\mu(g)=a $. 
Thus, if $\tau\in G$ is such that $a = \tau(0)$, and $\delta_{\tau}\in\mathcal M(G)$ denotes the Dirac mass at $\tau$, then $\mu_1 :=\check{\delta}_{\tau}*\mu*\delta_{\tau}$ satisfies $\int v(g)\;\text{d}\mu_1(g)=\int g(0)\;\text{d}\mu_1(g)=0$. Moreover, $\mu_1$ is compactly supported, and $\mu_1(\{g \, \vert \, v(g) \neq 0\}) > 0$, since $\mu$ satisfies (1).

Let $\ell\geq 1$  such that $\|\pi_0(\mu_1)\|_{L^2_0(\SO(d))}^{\ell}=\|\pi_0(\mu)\|_{L^2_0(\SO(d))}^{\ell}\leq 1/2$, and define $\mu_2 := \mu_1^{*(\ell)}$.
Then $\mu_2$ satisfies all the assumptions of Lemma \ref{small}. Therefore we can find a constant $c > 0$ such that for all $x \in \R^d$ satisfying $\vert x \vert \leq c$, we have
\[\Vert \pi_x(\mu_2) \Vert \leq 1 - c \vert x \vert^2.\]
Applying Lemma \ref{large} to some $y \in \R^d$ such that $\vert y \vert = c$, we find an absolute constant $C > 0$ such that for all $x \in \R^d$ satisfying $\vert x \vert \geq c$, we have
\[\Vert \pi_x(\mu_2)\Vert \leq 1 - C^{-1}c^3.\] 
Combining the last two displayed equalities yields a constant $c_1 > 0$ such that
\[\Vert \pi_x(\mu_2)\Vert \leq 1 - c_1 \min(\vert x \vert^2,1), \text{ for all } x \in \R^d.\]
Since $\pi_x(\mu_2) = \pi_x(\mu_1)^{\ell}=\pi_x(\tau)\pi_x(\mu)^{\ell}\pi_x(\tau)^{-1}$, we conclude that 
\[\|\pi_x(\mu)\| = \|\pi_x(\mu_2)\|^{1/\ell}\leq (1-c_1 \min(|x|^2,1))^{1/\ell} \leq 1 - \frac{c_1}{\ell}\min(\vert x \vert^2,1), \text{ for every } x\in\R^d.\]
This proves the result. \hfill$\blacksquare$
  
\section{Proofs of the main results}
In this section, we explain how to deduce Theorem \ref{B} from Theorem \ref{maintech}, and Theorem \ref{A} from Theorem \ref{B}. Finally, we prove Corollary \ref{C}.

\subsection{Proof of Theorem \ref{B}}  
Let $\Gamma < G$ be as in the statement of Theorem \ref{A}. We claim that there exists a symmetric finitely supported probability measure $\nu$ on $\Gamma$ satisfying the assumptions of Theorem \ref{maintech}.
Let $\mu \in \cM(G)$ be a symmetric measure whose support generates $\Gamma$. Then $\mu$ satisfies the assumptions of Theorem \ref{maintech}.
Let $\{F_n\}_n$ be an increasing sequence of symmetric finite subsets of $\Gamma$ such that $\cup_nF_n=\Gamma$.  For every $n$ with $\mu(F_n)>0$, define $\mu_n\in\cM(G)$ by letting $\mu_n(A)=\mu(A\cap F_n)/\mu(F_n)$, for any  set $A\subset G$. 
Then repeating the beginning of the proof of Theorem \ref{maintech} shows that $\nu:=\mu_n$ satisfies the assumptions of Theorem \ref{maintech}, for all large enough $n$.

Thus, we can find a constant $c_0 > 0$ such that $\Vert \pi_x(\nu) \Vert \leq 1 - c_0\min(\vert x \vert^2,1)$, for every $ x \in \R^d$. By assumption, the action $\theta(\Gamma) \actson \SO(d)$ has spectral gap, so we may take a convex combination of $\nu$ with another finitely supported measure on $\Gamma$ if necessary to assume that moreover, $\Vert \pi_0(\nu)\Vert_{L^2_0(\SO(d))} < 1$.

For every $\varphi\in L^2(\SO(d))$ and $x \in \R^d$, we have the inequality
\begin{equation}\label{est1}
\int_G\|\pi_x(g)\varphi-\varphi\|_2^2\;\text{d}\nu(g)=\langle (2-2\pi_x(\nu))(\varphi),\varphi\rangle\geq 2c_0\min\{|x|^2,1\}\|\varphi\|_2^2,
\end{equation}

Let $K\subset G$ be a compact set and $F \subset \Gamma$ denote the support of $\nu$. Let $K_0\subset\mathbb R^d$ be a compact set such that $K\cup F\subset K_0\cdot \SO(d)$.
Let $\kappa=\max\{|v|\;|\;v\in K_0\}$.
Note that $|2-e(y)-e(-y)|\leq 4\pi^2|y|^2$, for all $y\in\mathbb R$. 
Thus, for every $v\in K_0$ and $x\in\mathbb R^d$ we have that
\begin{align*}
\|\pi_x(v)\varphi-\varphi\|_2^2 &= \langle (2-\pi_x(v)-\pi_x(-v))(\varphi),\varphi\rangle\\
& =\int_{\SO(d)}\left(2-e(\langle \omega x,v\rangle) - e(-\langle \omega x,v\rangle)\right)\varphi(\omega)\overline{\varphi(\omega)}\;\text{d}\omega\\
&\leq 4\pi^2\kappa^2|x|^2\|\varphi\|_2^2.
\end{align*} 
Since $\|\pi_x(v)\varphi-\varphi\|_2^2\leq 4\|\varphi\|_2^2$, if we put $c_1=\max\{4\pi^2\kappa^2,4\}$, then the last inequality implies that 
\begin{equation}\label{est2}\|\pi_x(v)\varphi-\varphi\|_2^2\leq c_1\;\min\{|x|^2,1\}\|\varphi\|_2^2, \text{ for all } v\in K_0 \text{ and } x\in\mathbb R^d.\end{equation}
Denoting $c_2=c_1/2c_0$ and combining \eqref{est1} and \eqref{est2}, we conclude that 
\[\|\pi_x(v)\varphi-\varphi\|_2^2\leq c_2\int_G\|\pi_x(g)\varphi-\varphi\|_2^2\;\text{d}\nu(g), \text{ for all } v\in K_0 \text{ and } x\in\mathbb R^d.\]

Denote $\delta=\max_{g\in F}\|\lambda(g)\varphi-\varphi\|_2$.
Using that $\mathscr F\lambda\mathscr F^{-1}=\int_{\mathbb R^d}^{\oplus}\pi_x\;\text{d}x$ by equation \eqref{disint}, and integrating over $x \in \R^d$, we get
\begin{equation}\label{est3}\|\lambda(v)\varphi-\varphi\|_2^2\leq c_2\int_G\|\lambda(g)\varphi-\varphi\|_2^2\;\text{d}\nu(g)\leq c_2\delta^2, \text{ for all } v\in K_0, \text{ and } \varphi\in L^2(G).
\end{equation}

If $g\in F$, then $v(g)\in K_0$, hence $\|\lambda(\theta(g))(\varphi)-\varphi\|_2\leq\|\lambda(g)(\varphi)-\varphi\|_2+\|\lambda(v(g))\varphi-\varphi\|_2\leq (1+\sqrt{c_2})\delta$.
On the other hand, since $\|\pi_0(\nu)\|_{L^2_0(\SO(d))}<1$ and the restriction of $\lambda$ to $\SO(d)$ is a multiple of $\pi_0$, we can find a constant $c_3>0$, independent of $\varphi$, such that 
$$\sup_{\theta\in \SO(d)}\|\lambda(\theta)\varphi-\varphi\|_2\leq c_3\max_{g\in F}\|\lambda(\theta(g))\varphi-\varphi\|_2\leq c_3(1+\sqrt{c_2})\delta.$$

Thus, for every $g=(v,\theta)$, with $v\in K_0$ and $\theta\in \SO(d)$, we have that
$$\|\lambda(g)\varphi-\varphi\|_2\leq\|\lambda(v)\varphi-\varphi\|_2+\|\lambda(\theta)\varphi-\varphi\|_2\leq (\sqrt{c_2}+c_3+c_3\sqrt{c_2})\delta.$$

Since $K\subset K_0\cdot \SO(d)$, this proves the conclusion of Theorem B.\hfill$\blacksquare$

\subsection{Deducing Theorem \ref{A}}

The following lemma, inspired by the proof of \cite[Theorem 3]{Ma82}, explains why Theorem \ref{B} implies Theorem \ref{A}.  

\begin{lemma}\label{BtoA}
Let $\Gamma$ be a countable dense subgroup of an amenable locally compact second countable group $G$. Assume that for all sequences $\varphi_n\in L^2(G)$ satisfying $\|\lambda(g)\varphi_n-\varphi_n\|_2\rightarrow 0$, for all $g\in\Gamma$, and $\|\varphi_n\|_2=1$, for all $n$, we have that $\sup_{g\in K}\|\lambda(g)\varphi_n-\varphi_n\|_2\rightarrow 0$, for every compact set $K\subset G$.

Then the left translation action $\Gamma\curvearrowright G$ has local spectral gap.
\end{lemma}

{\it Proof.} Let $m$ be a left Haar measure of $G$, and $B\subset G$ be a measurable set with non-empty interior and compact closure.  In order to prove the conclusion, it suffices to argue that there is no sequence  $\xi_n\in L^2(G)$ such that $\|\xi_n\|_{2,B}=1$ and $\int_B\xi_n\;\text{d}m=0$, for every $n$, and $\|\lambda(g)\xi_n-\xi_n\|_{2,B}\rightarrow 0$, for every $g\in\Gamma$. 
Towards a contradiction, we prove that the ``mass" of $\xi_n$ becomes equidistributed:

{\bf Claim.} There is $\alpha>0$ such that
$\lim\limits_{n\rightarrow\infty}\int_GF|\xi_n|^2\;\text{d}m=\alpha\int_GF\;\text{d}m$, for every $F\in C_{\text{c}}(G)$.

{\it Proof of the claim.} It is clearly sufficient to show that the claim holds for a subsequence of $\{\xi_n\}$. 

Since $\|\lambda(g^{-1})\xi_n-\xi_n\|_{2,B}\rightarrow 0$ and $\|\xi_n\|_{2,B}=1$, we get that $\|\xi_n\|_{2,gB}\rightarrow 1$, thus $\sup_n\|\xi_n\|_{2,gB}<\infty$, for every $g\in\Gamma$.
If $C\subset G$ is a compact set, then $C$ can be covered with finitely many of the sets $\{gB\}_{g\in\Gamma}$, and therefore $\sup_n\|\xi_n\|_{2,C}<\infty$. This implies that we can find a Radon measure $\mu$ on $G$ and a subsequence $\{\xi_{n_k}\}$ of $\{\xi_n\}$ such that $\int_GF|\xi_{n_k}|^2\;\text{d}m\rightarrow\int_GF\;\text{d}\mu$, for every $F\in C_c(G)$.

Let $F\in C_{\text{c}}(G)$.  Denote by $C$ the support of $F$ and let $g\in\Gamma$. Since $C$ and $g^{-1}C$ can be covered with finitely many of the sets $\{hB\}_{h\in\Gamma}$, we have that $\|\lambda(g)\xi_n-\xi_n\|_{2,C}\rightarrow 0$. Since \begin{align*} \int_G(\lambda(g^{-1})F-F)|\xi_n|^2\;\text{d}m &=\int_GF\cdot(|\lambda(g)\xi_n|^2-|\xi_n|^2)\;\text{d}m\\ 
&\leq \|F\|_{\infty}\|\;|\lambda(g)\xi_n|^2-|\xi_n|^2\|_{1,C}\\&\leq \|F\|_{\infty}\|\lambda(g)\xi_n-\xi_n\|_{2,C}(\|\lambda(g)\xi_n\|_{2,C}+\|\xi_n\|_{2,C}) \end{align*} we deduce that $\int_G\lambda(g^{-1})F\;\text{d}\mu=\int_GF\;\text{d}\mu$. Since $F\in C_{\text{c}}(G)$ and $g\in\Gamma$ are arbitrary, we conclude that $\mu$ is $\Gamma$-invariant. Since $\Gamma<G$ is dense, $\mu$ must be $G$-invariant, and therefore there is $\alpha>0$ such that $\mu=\alpha\;m$. This finishes the proof of the claim. \hfill$\square$

Next, we claim that $1_C\xi_n\rightarrow 0$, weakly in $L^2(G)$, for every compact set $C\subset G$.  Indeed, since $\sup_n\|\xi_n\|_{2,C}<\infty$, for every $C\subset G$ compact, after replacing $\{\xi_n\}$ with a subsequence, we can find $\xi\in L^2_{\text{loc}}(G)$ such that $1_C\xi_n\rightarrow 1_C\xi$, weakly, for every $C\subset G$ compact. But then $\xi$ is $\Gamma$-invariant and since $\Gamma<G$ is dense, $\xi$ must be constant. On the other hand,  $\int_B\xi\;\text{d}m=\lim\limits_{n\rightarrow\infty}\int_B\xi_n\;\text{d}m=0$. By combining these facts we get that $\xi= 0$, almost everywhere, which proves our claim. 

Now, let $K\subset G$ be a non-negligible compact set. The hypothesis implies the existence of a finite set $S\subset\Gamma$ and $\delta>0$ such that if $\varphi\in L^2(G)$ satisfies $\max_{g\in S}\|\lambda(g)\varphi-\varphi\|_2\leq\delta\|\varphi\|_2$, then $\sup_{g\in K}\|\lambda(g)\varphi-\varphi\|_2\leq\|\varphi\|_2/2$. 

Finally, since $G$ is amenable we can find $F\in C_{\text{c}}(G)$ such that $\max_{g\in S}\|\lambda(g)F-F\|_2<\delta\|F\|_2$. Let $C$ be the support of $F$. Define $\varphi_n=F\xi_n\in L^2(G)$, for all $n$.  Then for every $g\in S$ we have \begin{align*} \|\lambda(g)\varphi_n-\varphi_n\|_2&\leq \|\lambda(g)F\cdot(\lambda(g)\xi_n-\xi_n)\|_2+\|(\lambda(g)F-F)\cdot\xi_n\|_2\\&\leq \|F\|_{\infty}\|\lambda(g)\xi_n-\xi_n\|_{2,gC}+\|(\lambda(g)F-F)\cdot\xi_n\|_2.\end{align*}

On the other hand, the claim implies that $\lim\limits_{n\rightarrow\infty}\|f\xi_n\|_2=\sqrt{\alpha}\|f\|_2$, for every $f\in C_{\text{c}}(G)$. In combination with the last inequality we conclude that for every $g\in S$ we have $$\limsup_{n\rightarrow\infty}\|\lambda(g)\varphi_n-\varphi_n\|_2\leq\sqrt{\alpha}\|\lambda(g)F-F\|_2<\delta\sqrt{\alpha}\|F\|_2=\delta\lim\limits_{n\rightarrow\infty}\|\varphi_n\|_2.$$
By the above, we derive that there is $N\geq 1$ such that $\sup_{g\in K}\|\lambda(g)\varphi_n-\varphi_n\|_2<\|\varphi_n\|_2/2$, for all $n\geq N$. Let $T:L^2(G)\rightarrow L^2(G)$ be the operator given by $\displaystyle{T(\varphi)=\frac{1_K}{m(K)}*\varphi}$.  Then 
\begin{equation}\label{T}\text{$\|T(\varphi_n)-\varphi_n\|_2<\|\varphi_n\|_2/2$, for all $n\geq N$.}\end{equation}
On the other hand, the restriction of $T$ to $L^2(C)$ is compact. Since $\varphi_n\in L^2(C)$, for every $n$, and $\varphi_n\rightarrow 0$, weakly, we deduce that $\|T(\varphi_n)\|_2\rightarrow 0$. In combination with \eqref{T} this implies that $\|\varphi_n\|_2\rightarrow 0$, which contradicts that $\|\varphi_n\|_2\rightarrow \sqrt{\alpha}\|F\|_2>0$.
\hfill$\blacksquare$

Although we will not need this later in the paper, we note that the converse of Lemma \ref{BtoA} also holds. We thank the referee for asking us whether this is the case.

\begin{lemma}\label{converse}
Let $\Gamma$ be a countable dense subgroup of a locally compact second countable group $G$. Assume that the left translation action $\Gamma\curvearrowright G$ has local spectral gap.

Then for every compact set $K\subset G$, we can find a finite set $F\subset\Gamma$ and a constant $\eta>0$ such that $\sup_{g\in K}\|\lambda(g)\varphi-\varphi\|_2\leq \kappa\max_{g\in F}\|\lambda(g)\varphi-\varphi\|_2$, for every $\varphi\in L^2(G)$.
\end{lemma}

Although Lemma \ref{converse} holds for arbitrary locally compact groups $G$, its conclusion is non-trivial only for $G$ amenable. Indeed, if $G$ is not amenable, then the restriction of $\lambda$ to any dense subgroup of $G$, and hence to $\Gamma$, has spectral gap. Thus, we can find a finite set  $F\subset\Gamma$  and $\eta>0$ such that  $\|\varphi\|_2\leq\eta\max_{g\in F}\|\lambda(g)\varphi-\varphi\|_2$, for all $\varphi\in L^2(G)$. This implies the conclusion of Lemma \ref{converse}.

{\it Proof.} Let $K\subset G$ be a compact set, and $m$ be a left Haar measure of $G$. Let $B_0\subset G$ be a compact set with non-empty interior. Since $B:=K^{-1}B_0\cup B_0$ is a compact set with non-empty interior, the action $\Gamma\curvearrowright G$ has local spectral gap with respect to $B$ (see \cite[Proposition 2.3]{BISG15}). By \cite[Proposition 2.2]{BISG15}, we can find 
a finite set $F\subset\Gamma$ and $\eta>0$ such that 
\begin{equation}\label{localsp}
\|\varphi-\frac{1}{m(B)}\int_B\varphi\;\text{d}m\|_{2,B}^2\leq\eta\sum_{h\in F}\|\lambda(h)\varphi-\varphi\|_{2,B}^2,\;\text{for every $\varphi\in L^2(G)$.}
\end{equation}
Fix $\varphi\in L^2(G)$ and denote $c=\frac{1}{m(B)}\int_B\varphi\;\text{d}m$. If $g\in K$, since $B_0\subset B$ and $g^{-1}B_0\subset B$, we get that \begin{align*}\|\lambda(g)\varphi-\varphi\|_{2,B_0}&\leq \|\lambda(g)\varphi-c\|_{2,B_0}+\|\varphi-c\|_{2,B_0}\\&=\|\varphi-c\|_{2,g^{-1}B_0}+\|\varphi-c\|_{2,B_0}\\&\leq 2\|\varphi-c\|_{2,B}.\end{align*}
By combining the last inequality and \eqref{localsp}, we deduce that 
\begin{equation}\label{localsp2}
\|\lambda(g)\varphi-\varphi\|_{2,B_0}^2\leq 4\eta\sum_{h\in F}\|\lambda(h)\varphi-\varphi\|_{2,B}^2,\;\text{for every $\varphi\in L^2(G)$ and $g\in K$.}
\end{equation}
Now, denote by $\rho:G\rightarrow\mathcal U(L^2(G))$ the right regular representation of $G$.
If $\varphi\in L^2(G)$, by using that $\lambda$ and $\rho$ commute and applying \eqref{localsp2} to $\rho(l)\varphi$, we get that 
\begin{equation}\label{localsp3}
\|\rho(l)(\lambda(g)\varphi-\varphi)\|_{2,B_0}^2\leq 4\eta\sum_{h\in F}\|\rho(l)(\lambda(h)\varphi-\varphi)\|_{2,B}^2,\;\text{for every $g\in K$ and $l\in G$.}
\end{equation}
Finally, $\int_G \|\rho(l)\psi\|_{2,C}^2\text{d}m(l)=\int_G\int_C|\psi(xl)|^2\text{d}m(x)\text{d}m(l)=\int_C\int_G|\psi(xl)|^2\text{d}m(l)\text{d}m(x)=m(C)\|\psi\|_2^2$, for every measurable set $C\subset G$ and $\psi\in L^2(G)$. In combination with \eqref{localsp3}, this gives that $m(B_0)\|\lambda(g)\varphi-\varphi\|_2^2\leq 4\eta \;m(B)\sum_{h\in F}\|\lambda(h)\varphi-\varphi\|_2^2$, for every $g\in K$, which finishes the proof.
\hfill$\blacksquare$

\subsection{ Proof of Corollary \ref{C}} By Theorem \ref{A}, the left translation action $\Gamma\curvearrowright G$ has local spectral gap with respect to any bounded measurable set $A\subset G$ with non-empty interior. Items (1) and (2) then follow from \cite[Theorem 7.1]{BISG15}, while item (3) follows from \cite[Theorem 4.1]{Io14}.

In order to prove (4), we will use the recent work of Grabowski, M\'athe, and Pikhurko \cite{GMP16}. Let $A, B\subset G$ be two bounded measurable subsets with non-empty interior.
By \cite[Corollary 4.5]{GMP16} the local spectral gap property can be rephrased as an expansion property. In particular, any  bounded measurable subset with non-empty interior, hence in particular $A$, is a domain of expansion for the action $\Gamma \actson G$. Thus, \cite[Theorem 1.5]{GMP16} implies that $A$ and $B$ are measurably $\Gamma$-equidecomposable if and only if they are $\Gamma$-equidecomposable and have the same Haar measure. 

Hence, we are left to prove that $A$ and $B$ are $\Gamma$-equidecomposable. Note that since the translation action of $\theta(\Gamma)$ on $\SO(d)$ has spectral gap, the group $\Gamma$ is non-amenable. Thus, since $G$ is connected, the assumptions of \cite[Lemma 7.11]{BISG15} are satisfied. This implies the desired conclusion.

For the moreover part, one can argue exactly the same way, provided that the analogue of \cite[Lemma 7.11]{BISG15} for the action $\Gamma \actson \R^d$ holds true. This is addressed in the next lemma.
\hfill$\blacksquare$

\begin{lemma}\label{equideco}
Let $G=\mathbb R^d\rtimes\SO(d)$, for some $d\geq 3$, and $\Gamma<G$ be a countable dense subgroup. 

Then any sets $A,B\subset\mathbb R^d$ with compact closure and non-empty interior are $\Gamma$-equidecomposable. 
\end{lemma} 

{\it Proof.} We follow the proof \cite[Lemma 7.11]{BISG15}. 
For every $r\geq 0$, let $B_r=\{x\in\mathbb R^d|\;|x|\leq r\}$. Using the argument from the end of the proof of \cite[Lemma 7.11]{BISG15} and the terminology therein, in order to get the conclusion, it suffices to argue that $B_1$ is $\Gamma$-paradoxical. 
Assume by contradiction that $B_1$ is not $\Gamma$-paradoxical. Then by a theorem of Tarski we can find a $\Gamma$-invariant, finitely additive measure $\varphi:\mathcal P(\mathbb R^d)\rightarrow [0,\infty]$ such that $\varphi(B_1)=1$.

Since $G$ is topologically perfect, 
 a result of Breuillard and Gelander \cite[Theorem 1.2]{BrG02} implies the existence of a finitely generated non-abelian free subgroup $\Gamma_0<\Gamma$ which is still dense in $G$.
 Moreover, the proof of \cite[Theorem 1.2]{BrG02}
 shows that for some fixed $m\geq 2$, for any neighbourhood $U$ of  the identity of $G$, we can find $b_1,...,b_m\in\Gamma\cap U$ such that  the group $\langle b_1,...,b_m\rangle$ generated by $b_1,...,b_m$ is isomorphic to $\mathbb F_m$ and dense in $G$.
 Thus, for every $1\leq i\leq m$, we can find a sequence $\{b_{n,i}\}_{n\geq 1}\subset\Gamma$ such that $b_{n,i}\rightarrow 1_G$ and the group $\langle b_{n,1},...,b_{n,m}\rangle$ is isomorphic to $\mathbb F_m$ and dense in $G$, for every $n\geq 1$.

Next, define $X=\sqcup_{n\geq 1}X_n$, with $X_n=\mathbb R^d$, to be the disjoint union of infinitely many copies of $\mathbb R^d$. Let $c_1,...,c_m$ denote the free generators of $\mathbb F_m$. We define an action $\mathbb F_m\curvearrowright X$ by letting $c_i\cdot x=b_{n,i}(x)$, for every $1\leq i\leq m$, $n\geq 1$ and $x\in X_n$. 

Let $s>0$ be a point at which the function $r\mapsto\varphi(B_r)$ is continuous. Define a finitely additive measure $\Phi:\mathcal P(\mathbb R^d)\rightarrow [0,1]$ by  $\Phi(A)=\varphi(A\cap B_s)/\varphi(B_s)$. Let $\omega\in\beta\mathbb N\setminus\mathbb N$ be a free ultrafilter on $\mathbb N$ and denote by $\lim\limits_{n\rightarrow\omega}x_n$ the limit of a bounded sequence $\{x_n\}_{n\geq 1}\subset\mathbb C$ along $\omega$. In particular,  $\lim\limits_{n\rightarrow\omega}$ is a bounded linear functional on $\ell^{\infty}(\mathbb N)$ which extends the usual limit. Define $\Psi:\mathcal P(X)\rightarrow [0,1]$ by letting $\Psi(A)=\lim\limits_{n\rightarrow\omega}\Phi(A\cap X_n)$. 
   Then the argument from the proof of \cite[Lemma 7.11]{BISG15} shows that $\Psi$ is an $\mathbb F_m$-invariant, finitely additive measure with $\Psi(X)=1$. In other words, $X$ admits an $\mathbb F_m$-invariant mean.

This further implies that the action $\mathbb F_m\curvearrowright X$ has the F$\o$lner property:
 there exists a sequence of finite subsets $\{F_k\}_{k\geq 1}$ of $X$ satisfying \begin{equation}\label{almost}\text{$\frac{|g\cdot F_k\triangle F_k|}{|F_k|}\rightarrow 0$, for every $g\in\mathbb F_2$}.\end{equation}
 
 This implication is proved in exactly the same way one deduces the existence of F$\o$lner sequences from the existence of a translation invariant mean for amenable groups (see e.g. \cite[Appendix G]{BdHV08}).
 Moreover, note that we may assume that for every $k\geq 1$, there is $n(k)\geq 1$ such that $F_k\subset X_{n(k)}=\mathbb R^d$. 

Next, we show that the sets $F_k$ are not ``trapped" in proper affine subspaces of $\mathbb R^d$.

{\bf Claim 1.} For $k\geq 1$, let $W_k\subset X_{n(k)}=\mathbb R^d$ be a proper affine subspace. Then $\lim\limits_{k\rightarrow\omega}\frac{|F_k\cap W_k|}{|F_k|}<1$.

{\it Proof of Claim 1.} Assume by contradiction that there exists a sequence of proper affine subspaces $W_k\subset X_{n(k)}$ such that $\lim\limits_{k\rightarrow\omega}\frac{|F_k\cap W_k|}{|F_k|}=1$. Moreover, we can choose the sequence $\{W_k\}_{k\geq 1}$ such that $\rho:=\lim\limits_{k\rightarrow\omega}\dim(W_k)\in\{0,1,...,d-1\}$ is minimal among all such sequences.

Let $g\in\mathbb F_m$. Then $|F_k\cap g\cdot W_k|=|g^{-1}\cdot F_k\cap W_k|\geq |F_k\cap W_k|-|g^{-1}F_k\triangle F_k|$, and by using \ref{almost} we deduce that $\lim\limits_{k\rightarrow\omega}\frac{|F_k\cap g\cdot W_k|}{|F_k|}=1$. This further gives that 
 $\lim\limits_{k\rightarrow\omega}\frac{|F_k\cap (W_k\cap g\cdot W_k)|}{|F_k|}=1$. The minimality assumption on $\{W_k\}_{k\geq 1}$ implies that $\lim\limits_{k\rightarrow\omega}\dim(W_k\cap g\cdot W_k)=\rho$, for every $g\in\mathbb F_m$.
 
  By applying this fact to $g\in\{c_1,...,c_m\}$ we get $\lim\limits_{k\rightarrow\omega}\sum_{i=1}^m\big(\dim(W_k)-\dim(W_k\cap c_i\cdot W_k)\big)=0$. In particular, we can find $k\geq 1$ such that $\sum_{i=1}^m\big(\dim(W_k)-\dim(W_k\cap c_i\cdot W_k)\big)=0$. Hence $W_k=c_i\cdot W_k$, for every $1\leq i\leq m$. Thus, $W_k\subset X_{n(k)}=\mathbb R^d$ is a proper affine subspace which is invariant under the group of isometries $\langle b_{n(k),1},...,b_{n(k),m}\rangle.$
 Since the latter is dense in $G$, we deduce that $W_k$ is invariant under $G$, which is a contradiction. \hfill$\square$
 
For $k\geq 1$, let $S_k$ be the  set of $(x_1,...,x_{d+1})\in F_k^{d+1}$ such that $x_1,...,x_{d+1}$ are affinely independent.

{\bf Claim 2.} $\lim\limits_{k\rightarrow\omega}\frac{|S_k|}{|F_k^{d+1}|}>0$.

{\it Proof of Claim 2.} For $k\geq 1$ and $1\leq p\leq d+1$, we denote by $S_{k,p}$ the set of $(x_1,...,x_p)\in F_k^p$ such that $x_1,...,x_p$ are affinely independent. 
Let $1\leq p\leq d+1$ be the largest integer such that $\lim\limits_{k\rightarrow\omega}\frac{|S_{k,p}|}{|F_k^p|}>0$.
Our goal is to show that $p=d+1$.
Assume by contradiction that $p\leq d$. 

Then $\lim\limits_{k\rightarrow\omega}\frac{|S_{k,p+1}|}{|F_k^{p+1}|}=0$. Since $\lim\limits_{k\rightarrow\omega}\frac{|S_{k,p}|}{|F_k^p|}>0$, we get that $\lim\limits_{k\rightarrow\omega}\frac{|S_{k,p+1}|}{|S_{k,p}\times F_k|}=0$. Let $\varepsilon_k\in [0,1]$ such that  $|S_{k,p+1}|=\varepsilon_k |S_{k,p}\times F_k|$. Then $\lim\limits_{k\rightarrow\omega}\varepsilon_k=0$ and for every $k\geq 1$ we have $$\varepsilon_k|S_{k,p}|\;|F_k|=|S_{k,p+1}|=\sum_{(x_1,...,x_p)\in S_{k,p}}|\{x_{p+1}\in F_k\;|\;x_1,...,x_p,x_{p+1}\;\text{are affinely independent}\}|.$$

Let $(x_1,...,x_p)\in S_{k,p}$ such that 
$|\{x_{p+1}\in F_k\;|\;x_1,...,x_p,x_{p+1}\;\text{are affinely independent}\}|\leq\varepsilon_k|F_k|$.
Denote by $W_k\subset\mathbb R^d$  the affine subspace spanned by $x_1,...,x_p$. Then $|F_k\cap W_k|\geq (1-\varepsilon_k) |F_k|$ and thus $\lim\limits_{k\rightarrow\omega}\frac{|F_k\cap W_k|}{|F_k|}=1$. Since $p\leq d$, the subspace $W_k\subset\mathbb R^d$ is proper, which contradicts Claim 1. \hfill$\square$

For every $n\geq 1$, let $\tilde X_n:=\{(x_1,...,x_{d+1})\in X_n^{d+1}\;|\; x_1,...,x_{d+1}\;\text{are affinely independent}\}$. Put $\tilde X=\sqcup_{n\geq 1}\tilde X_n$, and consider the diagonal action $\mathbb F_m\curvearrowright \tilde X$ given by $g\cdot x=(g\cdot x_1,...,g\cdot x_{d+1})$, for every $g\in\mathbb F_m$, $n\geq 1$ and $x=(x_1,...,x_{d+1})\in\tilde X_n$. Note that the action $\mathbb F_m\curvearrowright \tilde X$ is free.

Recall that $F_k\subset X_{n(k)}$ and observe that $S_{k}=F_k^{d+1}\cap \tilde X_{n(k)}$. If $g\in\mathbb F_m$, then $g\cdot S_k=g\cdot F_k^{d+1}\cap \tilde X_{n(k)}$. Hence $g\cdot S_k\triangle S_k\subset g\cdot F_k^{d+1}\triangle F_k^{d+1}$ and thus $$\frac{|g\cdot S_k\triangle S_k|}{|S_k|}\leq \frac{|g\cdot F_k^{d+1}\triangle F_k^{d+1}|}{|F_k^{d+1}|}\;\;\frac{|F_k^{d+1}|}{|S_k|}\leq (d+1)\frac{|g\cdot F_k\triangle F_k|}{|F_k|}\;\;\frac{|F_k^{d+1}|}{|S_k|}$$. 

By combining \ref{almost} and Claim 2, we deduce that $\lim\limits_{k\rightarrow\omega}\frac{|g\cdot S_k\triangle S_k|}{|S_k|}=0$, for every $g\in\mathbb F_m$.
In particular, we can find a subsequence $\{S_{k(\ell)}\}_{\ell\geq 1}$ of $\{S_k\}_{k\geq 1}$ such that $\lim\limits_{\ell\rightarrow\infty}\frac{|g\cdot S_{k(\ell)}\triangle S_{k(\ell)}|}{|S_{k(\ell)}|}=0$, for every $g\in\mathbb F_m$.

However, since the action $\mathbb F_m\curvearrowright\tilde X$ is free, this contradicts the non-amenability of $\mathbb F_m$.
\hfill$\blacksquare$


\end{document}